\documentclass[11pt,fleqn]{article}
\usepackage{stmaryrd}
\usepackage{amssymb}
\usepackage{bbm}
\usepackage{mathrsfs}
\usepackage{amsfonts}
\usepackage{amsmath,epsf,eufrak}
\begin{document}
\renewcommand{\theequation}{\arabic {section}.\arabic {equation}}
\newcommand{\bi}{\begin{equation}}
\newcommand{\ei}{\end{equation}}
\date{}
\title{  On the Sendov conjecture and the critical points of polynomials
\author{Zaizhao Meng}}
\maketitle \baselineskip 24pt
 \begin{center}
 In this paper, we obtain new results on the critical points of a polynomial. We discuss the Sendov conjecture for polynomials of degree nine.
 \end{center}
\section{Introduction}
\setcounter{equation}{0}

Let $ \mathcal{P}_{n}$ denote the set of all monic polynomials of degree $n(\geq 2)$ of the form
$$ p(z)=\prod\limits_{k=1}^{n}(z-z_{k}),\ \ |z_{k}|\leq 1(k=1,\cdots,n)$$
with\\
$p^{\prime}(z)=n\prod\limits_{j=1}^{n-1}(z-\zeta_{j}),\ \ |\zeta_{j}|\leq 1(j=1,\cdots,n-1).$\\
Write $I(z_{k})=\min\limits_{1\leq j\leq n-1}|z_{k}-\zeta_{j}|,\ \ I(p)=\max\limits_{1\leq k\leq n}I(z_{k})$,
and $I(\mathcal{P}_{n})=\sup\limits_{p\in\mathcal{P}_{n}}I(p)$.\\
It was showed that there exists an extremal polynomial $p^{\ast}_{n}$, i.e., $I(\mathcal{P}_{n})=I(p^{\ast}_{n})$ and that $p^{\ast}_{n}$ has at least one zero on
each subarc of the unit circle of length $\pi$(see [3], [9]).\\
It will suffice to prove the Sendov conjecture assuming $p$ is an extremal polynomial of the following form,\\
$p(z)=(z-a)\prod\limits_{k=1}^{n-1}(z-z_{k}),\ \ |z_{k}|\leq 1(k=1,\cdots,n-1)$\\
with\\
$p^{\prime}(z)=n\prod\limits_{j=1}^{n-1}(z-\zeta_{j}),\ \ |\zeta_{j}|\leq 1(j=1,\cdots,n-1),\ \ a\in [0,1]$.\\
Let $r_{k}=|a-z_{k}|,\ \rho_{j}=|a-\zeta_{j}|$ for
$k,j=1,2,\cdots,n-1$. By relabeling we suppose that \bi
\rho_{1}\leq\rho_{2}\leq \cdots\leq\rho_{n-1},\  r_{1}\leq r_{2}\leq\cdots\leq r_{n-1} . \ei
 We have(see [5],[7],[10])
 \bi
 2\rho_{1}\sin(\frac{\pi}{n})\leq r_{k}\leq 1+a,\ \
k=1,2,\cdots,n-1.
 \ei
 {\bf Sendov conjecture.} The disk $|z-a|\leq 1$ contains a zero of $p^{\prime}(z)$.\\
 If $a=0$  or $a=1$, Sendov conjecture is true(see [7],[11]), we suppose $a\in (0,1)$.\\
 In this paper, we obtain the following theorems.\\
 {\bf Theorem 1.} If
$\ \ 1-(1-|p(0)|)^{\frac{1}{n}}\leq \lambda\leq\sin(\frac{\pi}{n})$ and
$\lambda<a,\ \rho_{1}\geq 1$, then there exists a critical point
$\zeta_{0}=a+\rho_{0}e^{i\theta_{0}}$ such that $ Re \zeta_{0}\geq
\frac{1}{2}(a-\frac{\lambda(\lambda+2)}{a})$.\\
This theorem improves the previous known results(see [2],[4],[5],[8]).\\
{\bf Theorem 2.} If $p$ is an extremal polynomial, $n=9,\ p(a)=0$, and $\ a\in [0,0.845)$, then the disk $|z-a|\leq 1$ contains a zero of $p^{\prime}(z)$.\\
{\bf Theorem 3.} For $\rho>0,$ and $ m $\ is real number, we have\\
$\prod\limits_{r_{k}\geq \rho}r_{k}\rho^{-1}\leq \prod\limits_{\rho_{j}\geq \rho^{m}}\rho_{j}\rho^{-m}
\prod\limits_{\rho^{m-1}2\sin\frac{\pi k}{n}\geq 1}\rho^{m-1}2\sin\frac{\pi k}{n}.$

\section{ Proof of the Theorem 1 }
\setcounter{equation}{0}
{\bf Lemma 2.1.} If $0<a<1$ and $\rho_{1}\geq 1$ , then\\
$|p(z)|>1-(1-\lambda)^{n}$ ,\ for
$0<|z-a|=\lambda\leq\sin(\pi/n)$.\\
This is Lemma 2.16 of [2](see [1]).\\
{\bf Proof of the Theorem 1.}\ \
We apply Lemma 2.1 to conclude that
\\ $|p(z)|>1-(1-\lambda)^{n}\geq|p(0)|,\ 0<|z-a|=\lambda$.\\
Since $p(z)$ is univalent in $|z-a|\leq\lambda$, it follows that
there exists a unique point $z_{0}$ with $|z_{0}-a|<\lambda$ such
that $p(0)=p(z_{0})$. We assume that $Im z_{0}\geq 0$(if not,
consider $\overline{p(\overline{z})}$). By a variant of the
Grace-Heawood theorem, there exists a critical point in each of the
half-planes bounded by the perpendicular bisector $L$ of the segment
from 0 to $z_{0}$. Let $ \zeta_{0}=a+\rho_{0}e^{i\theta_{0}}$ be the
critical point in the half-plane containing $z_{0}$.\\
The equation of $L$ is $|z|=|z-z_{0}|$, that is
 \bi
z\overline{z}=(z-z_{0})(\overline{z}-\overline{z_{0}}),
 \ei
 then
 \bi
 |z_{0}|^{2}=z \overline{z_{0}}+\overline{z}z_{0}.
\ei Let $z^{\ast}=e^{i\beta_{0}}$ be the joint point of $L$ and the
circle $|z|=1,\ Im z^{\ast}\geq 0$, then\\
$|z_{0}|^{2}=e^{i\beta_{0}}\overline{z_{0}}+e^{-i\beta_{0}}z_{0}$.\\
Hence\\
$e^{i\beta_{0}}=\frac{z_{0}}{2}\pm\frac{z_{0}}{2|z_{0}|}\sqrt{|z_{0}|^{2}-4}$,\\
that is
\bi
e^{i\beta_{0}}=(\frac{1}{2}\pm\frac{i}{2|z_{0}|}\sqrt{4-|z_{0}|^{2}})z_{0}.
\ei
If $z_{0}=a$, the theorem is true, we write $z_{0}=a+re^{i\alpha},\ 0<r<\lambda$. We choose\\
$\cos\beta_{0}=\frac{1}{2}(a+r\cos\alpha)-\frac{1}{2|z_{0}|}\sqrt{4-|z_{0}|^{2}}
r\sin\alpha, $
\bi
\cos\beta_{0}=\frac{1}{2}(a+r\cos\alpha)-\frac{1}{2}\frac{\sqrt{4-a^{2}-2ar\cos\alpha-r^{2}}}{\sqrt{a^{2}+2ar\cos\alpha+r^{2}}}r\sin\alpha.
\ei
We fix $r$ and consider the circle $|z-a|=r$, let
$\sin\alpha_{1}=\frac{\sqrt{a^{2}-r^{2}}}{a} , \ x=\cos\alpha$ and $
F(x)=x-\sqrt{\frac{4-a^{2}-2arx-r^{2}}{a^{2}+2arx+r^{2}}}\sqrt{1-x^{2}}$,\\
then
\bi
 \cos\beta_{0}=\frac{1}{2}(a+rF(x)),
\ei
it is sufficient to give lower bound of $F(x)$ for
$x\in[-1,-\frac{r}{a}]$.\\
Write $G(x)=-F(-x)$, then
\bi
G(x)=x+\sqrt{\frac{4-a^{2}+2arx-r^{2}}{a^{2}-2arx+r^{2}}}\sqrt{1-x^{2}},
\ei
it is sufficient to give upper bound of $G(x)$ for $x\in
[\frac{r}{a},1]$.\\
Write $\psi=\frac{4-a^{2}+2arx-r^{2}}{a^{2}-2arx+r^{2}},\
\phi=a^{2}-2arx+r^{2}.$\\
We have
\bi
G^{\prime}(x)=1+\frac{1}{2}\psi^{-\frac{1}{2}}\psi^{\prime}\sqrt{1-x^{2}}-\psi^{\frac{1}{2}}x(1-x^{2})^{-\frac{1}{2}},
\ei
and $\psi^{\prime}=8ar(a^{2}-2arx+r^{2})^{-2}$,\\
hence
\bi
 G^{\prime}(x)=\phi^{-2}\psi^{-\frac{1}{2}}(1-x^{2})^{-\frac{1}{2}}G_{1},
\ei
where
\bi
G_{1}=\phi^{\frac{3}{2}}(4-\phi)^{\frac{1}{2}}(1-x^{2})^{\frac{1}{2}}+4ar(1-x^{2})-x\phi(4-\phi).
\ei
We have $x=\frac{a^{2}+r^{2}-\phi}{2ar}$, then
\bi
 G_{1}=\frac{1}{2ar}G_{2},
\ei
where\\
$G_{2}=\phi^{\frac{3}{2}}(4-\phi)^{\frac{1}{2}}(4a^{2}r^{2}-(a^{2}+r^{2})^{2}+2(a^{2}+r^{2})\phi-\phi^{2})^{\frac{1}{2}}\\
+8a^{2}r^{2}-2(a^{2}+r^{2})^{2}+4(a^{2}+r^{2})\phi-2\phi^{2}-(a^{2}+r^{2}-\phi)\phi(4-\phi)$, \\
hence\\
$G_{2}=\phi^{\frac{3}{2}}(4-\phi)^{\frac{1}{2}}(4a^{2}r^{2}-(a^{2}+r^{2})^{2}+2(a^{2}+r^{2})\phi-\phi^{2})^{\frac{1}{2}}\\
+8a^{2}r^{2}-2(a^{2}+r^{2})^{2}+(a^{2}+r^{2}+2)\phi^{2}-\phi^{3},$\\
and $\phi\in [(a-r)^{2},a^{2}-r^{2}].$\\
The roots of $G_{2}$ satisfy\\
$L=\phi^{3}(4-\phi)(4a^{2}r^{2}-(a^{2}+r^{2})^{2}+2(a^{2}+r^{2})\phi-\phi^{2})=\\
(\phi^{3}-(a^{2}+r^{2}+2)\phi^{2}+2(a^{2}+r^{2})^{2}-8a^{2}r^{2})^{2}=R,$ \\
say.\\
We have\\
$L=\phi^{6}-2(a^{2}+r^{2}+2)\phi^{5}+(8(a^{2}+r^{2})-4a^{2}r^{2}+(a^{2}+r^{2})^{2})\phi^{4}+
(16a^{2}r^{2}-4(a^{2}+r^{2})^{2})\phi^{3}$,\\
$R=\phi^{6}-2(a^{2}+r^{2}+2)\phi^{5}+(a^{2}+r^{2}+2)^{2}\phi^{4}+4(a^{2}-r^{2})^{2}\phi^{3}
-4(a^{2}+r^{2}+2)(a^{2}-r^{2})^{2}\phi^{2}+4(a^{2}-r^{2})^{4}$.\\
Let $d=1-r^{2}, \ e_{1}=a^{2}-1$, by $L=R$, we deduce\\
$e_{1}d\phi^{4}-2(e_{1}+d)^{2}\phi^{3}+(4+e_{1}-d)(e_{1}+d)^{2}\phi^{2}-(e_{1}+d)^{4}=0$,\\
that is\\
$(e_{1}\phi^{2}-2(e_{1}+d)\phi-(e_{1}+d)^{2})(d\phi^{2}-2(e_{1}+d)\phi+(e_{1}+d)^{2})=0$,\\
there is only the root $\phi_{0}=\frac{a^{2}-r^{2}}{1+r}\in [(a-r)^{2},a^{2}-r^{2}]$, and\\
$x_{0}=\frac{a^{2}+r^{2}-\phi_{0}}{2ar}=\frac{2r+a^{2}+r^{2}}{2a(1+r)}\in
[\frac{r}{a},1]$.\\
We have $G^{\prime}(1-0)<0,\ G^{\prime}(\frac{r}{a})>0$, $G(x_{0})$ is
the maxima value of $G(x)$ for $x\in [\frac{r}{a},1]$.\\
We obtain  $G(x_{0})=\frac{r+2}{a}$,and $G(x)\leq\frac{r+2}{a}, \ x\in [\frac{r}{a},1]$, hence\\
 $F(x)\geq -\frac{r+2}{a},\ x\in [-1,-\frac{r}{a}],\\
\cos\beta_{0}\geq\frac{1}{2}(a-\frac{r(r+2)}{a})\geq\frac{1}{2}(a-\frac{\lambda(\lambda+2)}{a}),$\\
by the Gauss-Lucas theorem , Theorem 1 follows.
\section{ Proof of the Theorem 2 }
\setcounter{equation}{0}
We improve the methods of [5].\\
If $p(0)=0$, the Sendov conjecture is true(see Satz 3 of [11]).\\
In this section, we assume that $z_{k}\neq 0, \ \ z_{k}\neq a,\ k=1,\cdots,n-1$ and that $p(z)$ is extremal: $I(\mathcal{P}_{n})=I(p)=I(a)=\rho_{1}.$\\
{\bf Lemma 3.1.} If $c_{k}(k=1,\cdots,N),\ m,\ M, \ C$ are positive constants with\\
$m\leq c_{k}\leq M,\ \prod\limits_{k=1}^{N}c_{k}\geq C$ and $m^{N}\leq C\leq M^{N}$, then
$$\sum\limits_{k=1}^{N}\frac{1}{c_{k}^{2}}\leq\frac{N-v}{m^{2}}+\frac{v-1}{M^{2}}+\{\frac{m^{N-v}M^{v-1}}{C}\}^{2},$$
where $v=\min\{j\in \mathbb{Z}:M^{j}m^{N-j}\geq C\}.$\\
Proof. See Lemma 7.3.9 of [10].\\
{\bf Lemma 3.2.} If $w\neq a$, then
$$ \prod\limits_{k=1}^{n-1}(a-z_{k})=n\prod\limits_{j=1}^{n-1}(a-\zeta_{j}),$$
$$ \prod\limits_{k=1}^{n-1}r_{k}=n\prod\limits_{j=1}^{n-1}\rho_{j},$$
$$\sum\limits_{j=1}^{n-1} \frac{1}{a-\zeta_{j}}= \sum\limits_{k=1}^{n-1} \frac{2}{a-z_{k}},$$
$$Re\{\frac{1}{a-w}\}=\frac{1}{2a}-\frac{|w|^{2}-a^{2}}{2a|a-w|^{2}},$$
$$ \frac{1}{n}(a+\sum\limits_{k=1}^{n-1}z_{k})=\frac{1}{n-1}\sum\limits_{j=1}^{n-1}\zeta_{j}.$$
Proof. See (2.2), (2.3), (2.5) of [5], and (2.3.1) of [10] .\\
Let $$\gamma_{j}=\frac{\zeta_{j}-a}{a\zeta_{j}-1}\  and \ w_{k}=\frac{z_{k}-a}{az_{k}-1}. $$
By (2.11) of [5], we have
$$\prod\limits_{j=1}^{n-1}|\gamma_{j}|\leq \frac{\prod\limits_{k=1}^{n-1}|w_{k}|}{n-\frac{4a^{2}}{1+a^{2}}-a\sum\limits_{k=1}^{n-3}Re w_{k}}.$$
We take $n=9$, then
\bi
\prod\limits_{j=1}^{8}|\gamma_{j}|\leq \frac{\prod\limits_{k=1}^{8}|w_{k}|}{9-4a^{2}/(1+a^{2})-6a}.
\ei
{\bf Lemma 3.3.} Let $A_{9}(A_{9}=0.4314\cdots)$ be the smallest positive root of $9-\frac{4x^{2}}{1+x^{2}}-6x-(1+x-x^{2})^{8}=0$, and
$a\leq A_{9}$, then $\rho_{1}\leq 1$.\\
Proof. This is Lemma 3.2 of [5].\\
{\bf Lemma 3.4.} If $\rho_{1}> 1,\ $ and $\ \zeta_{0}=a+\rho_{0}e^{i\theta_{0}}$  is the critical point in Theorem 1,  $\gamma_{0}=\frac{\zeta_{0}-a}{a\zeta_{0}-1}$, then
$$ |\gamma_{0}|>\frac{1}{\sqrt{1+\lambda(\lambda+2)-a^{2}\lambda(\lambda+2)}}.$$
Proof. We have $ |\gamma_{0}|^{2}=\frac{\rho_{0}^{2}}{a^{2}|\zeta_{0}|^{2}-2aRe\zeta_{0}+1},\ Re \zeta_{0}\geq
\frac{1}{2}(a-\frac{\lambda(\lambda+2)}{a})$, hence\\
$\rho_{0}\cos\theta_{0}\geq-\frac{1}{2}(a+\frac{\lambda(\lambda+2)}{a})$, and \\
$ |\gamma_{0}|^{2}=\frac{\rho_{0}^{2}}{a^{4}-2a^{2}+1+a^{2}\rho_{0}^{2}+2a(a^{2}-1)\rho_{0}\cos\theta_{0}}
\geq\frac{\rho_{0}^{2}}{a^{4}-2a^{2}+1+a^{2}\rho_{0}^{2}+(1-a^{2})(a^{2}+\lambda(\lambda+2))}\\
>\frac{1}{1+\lambda(\lambda+2)-a^{2}\lambda(\lambda+2)},$\\
the lemma follows.\\
{\bf Lemma 3.5.} If $|\gamma_{j}|\leq\frac{1}{1+a-a^{2}}$, then $\rho_{j}\leq 1$.\\
Proof. See Lemma 1 of [4].\\
{\bf Lemma 3.6.} If $R\in (0,1],\  w=\frac{z-a}{az-1},\ |z|\leq R$, then
$$|w|\leq \frac{|z|+a}{a|z|+1}\leq\frac{R+a}{aR+1}.$$
Proof. We have \\
$|w|^{2}=\frac{|z-a|^{2}}{|az-1|^{2}}=\frac{|z|^{2}-2aRe z+a^{2}}{a^{2}|z|^{2}-2aRe z+1}\leq\frac{|z|^{2}+2a|z|+a^{2}}{a^{2}|z|^{2}+2a|z|+1}=\frac{(|z|+a)^{2}}{(a|z|+1)^{2}},$\\
 hence\\
$|w|\leq \frac{|z|+a}{a|z|+1}\leq\frac{R+a}{aR+1}$.\\
Write $B=B(R)=\frac{R+a}{aR+1}$.\\
{\bf Lemma 3.7.} If $ a\in [0.4314,0.51952]$, or if there exists $|z_{k}|\leq 0.4$ for $ a\in [0.5195,1]$,
 then $I(a)=\rho_{1}\leq 1$.\\
Proof. If $I(a)>1$ and there exists $|z_{k}|\leq R$, then,
by (3.1), Lemma 3.4, and Lemma 3.5, there exists some $\gamma_{j_{0}}$,\\
$ \frac{|\gamma_{j_{0}}|^{7}}{\sqrt{1+\lambda(\lambda+2)-a^{2}\lambda(\lambda+2)}}<\prod\limits_{j=1}^{8}|\gamma_{j}|< \frac{B}{9-4a^{2}/(1+a^{2})-6a},$\\
hence\\
$|\gamma_{j_{0}}|^{7}<\frac{\sqrt{1+\lambda(\lambda+2)-a^{2}\lambda(\lambda+2)}}{9-4a^{2}/(1+a^{2})-6a}B$.\\
By Lemma 3.5, it suffices to show
\bi
 (9-4a^{2}/(1+a^{2})-6a)\frac{aR+1}{ R+a}- \sqrt{1+(1-a^{2})\lambda(\lambda+2)}(1+a-a^{2})^{7}\geq 0.
\ei
We consider the conditions for $\lambda$,\\
$|p(0)|=a\prod\limits_{k=1}^{8}|z_{k}|\leq aR,$
$\ \ 1-(1-|p(0)|)^{\frac{1}{9}}\leq \lambda\leq\sin(\frac{\pi}{9}),$\\
we choose
$$\lambda=1-(1-aR)^{\frac{1}{9}},$$
and $R$ satisfies $R\leq a^{-1}(1-(1-\sin(\frac{\pi}{9}))^{9}).$\\
If $ a\in [0.4314,0.51952]$, we take $R=1$, then
$$9-4a^{2}/(1+a^{2})-6a-(1+a-a^{2})^{7}\sqrt{1+(1-a^{2})\lambda(\lambda+2)}>0,$$
we obtain (3.2).\\
If $ a\in [0.5195,1]$, and there exists $|z_{k}|\leq 0.4$, we take $R=0.4$, (3.2) holds again,
the lemma follows.\\
{\bf Lemma 3.8.}  If $a\in(0.5195,1],\ x\in[0.4,1]$, then \\
$x^{2}+\frac{(1-x^{2})(a^{2}+x^{2})^{2}}{9x^{2}}\leq 1$.\\
Proof. It is sufficient to show
$$ (1+x^{2})^{2}\leq 9x^{2},$$
this holds for $x\in[0.4,1]$, the lemma follows.\\
We have $\frac{p^{\prime}(0)}{p(0)}=-(\frac{1}{a}+\sum\limits_{k=1}^{8}\frac{1}{z_{k}}),$  and
$9\prod\limits_{j=1}^{8}|\zeta_{j}|=(a\prod\limits_{k=1}^{8}|z_{k}|)|\frac{1}{a}+\sum\limits_{k=1}^{8}\frac{1}{z_{k}}|.$\\
let $\Delta=Re(\frac{1}{a}+\sum\limits_{k=1}^{8}\frac{1}{z_{k}}),
\ \sigma=\sum\limits_{k=1}^{8}\frac{1}{r_{k}^{2}}$,
then $9\prod\limits_{j=1}^{8}|\zeta_{j}|\geq -\Delta a\prod\limits_{k=1}^{8}|z_{k}|$.\\
{\bf Lemma 3.9.} If $\rho_{1}>1,\ a\in [0.5195,1]$, and $|z_{k}|\in [0.4, 1], \ k=1,\cdots,8$,  then
$$\Delta\leq -\frac{8}{a}+8a+\frac{9}{8a}(1-a^{2})\sigma.$$
Proof. Write $z_{k}=|z_{k}|e^{i\theta_{k}},\ \zeta_{j}=a+\rho_{j}e^{it_{j}},\  1\leq k,\ j\leq 8.$
By Lemma 3.2, we have\\
$ \Delta=\frac{1}{a}+Re\sum\limits_{k=1}^{8}z_{k}+\sum\limits_{k=1}^{8}\frac{1-|z_{k}|^{2}}{|z_{k}|}\cos\theta_{k},$
\bi
\Delta=\frac{1-a^{2}}{a}+\frac{9}{8}Re\sum\limits_{j=1}^{8}\zeta_{j}+
\sum\limits_{k=1}^{8}\frac{1-|z_{k}|^{2}}{|z_{k}|}\cos\theta_{k},
\ei
$$\sum\limits_{k=1}^{8}\frac{|z_{k}|^{2}-a^{2}}{r_{k}^{2}}=8+a\sum\limits_{j=1}^{8}\frac{\cos t_{j}}{\rho_{j}}. $$
Since $\rho_{j}>1$, we deduce $\cos t_{j}\leq 0,$\\
$\sum\limits_{k=1}^{8}\frac{|z_{k}|^{2}-a^{2}}{r_{k}^{2}}\geq 8+a\sum\limits_{j=1}^{8}\cos t_{j}. $\\
But \\
$Re\sum\limits_{j=1}^{8}\zeta_{j}=8a+\sum\limits_{j=1}^{8}\rho_{j}\cos t_{j}\leq 8a+\sum\limits_{j=1}^{8}\cos t_{j}$,\\
hence\\
$Re\sum\limits_{j=1}^{8}\zeta_{j}\leq 8a-\frac{8}{a}+\frac{1}{a}\sum\limits_{k=1}^{8}\frac{|z_{k}|^{2}-a^{2}}{r_{k}^{2}}$,
and by (3.3)
\bi
\Delta\leq 8a-\frac{8}{a}+\sum\limits_{k=1}^{8}(\frac{9}{8a}\frac{|z_{k}|^{2}-a^{2}}{r_{k}^{2}}
+\frac{1-|z_{k}|^{2}}{|z_{k}|}\cos\theta_{k}).
\ei
We want to prove
\bi
\frac{9}{8a}\frac{|z_{k}|^{2}-a^{2}}{r_{k}^{2}}
+\frac{1-|z_{k}|^{2}}{|z_{k}|}\cos\theta_{k}\leq \frac{9}{8a}\frac{1-a^{2}}{r_{k}^{2}}.
\ei
If $\cos\theta_{k}\leq 0$, (3.5) is valid , we assume $\cos\theta_{k}>0$ and write $x=|z_{k}|,\  \theta=\theta_{k}\ $,
then  $r_{k}^{2}=a^{2}+x^{2}-2ax\cos\theta$. We want to show

$$x^{2}+\frac{8a (1-x^{2})}{9 x}(a^{2}+x^{2}-2ax\cos\theta)\cos\theta\leq 1,  $$
it suffices to prove
$$
x^{2}+\frac{(1-x^{2})(a^{2}+x^{2})^{2}}{9 x^{2}}\leq 1,
$$
for $ a\in(0.5195,1],\ x\in[0.4,1]$, this is true by Lemma 3.8, the lemma follows.\\
{\bf Lemma 3.10.} Let $m=\frac{1}{4},\ a\in [0.5,0.95], \ f(x)=\frac{x^{2}-1}{(1-x^{m})(a+x)^{2}}$, then
$f^{\prime}(x)>0$, for $x\in [0.4,1)$.\\
Proof. We have  $f^{\prime}(x)=\frac{a+x}{(1-x^{m})^{2}(a+x)^{4}}Y$, where \\
$Y=((m-2)x^{m+1}-mx^{m-1}+2x)a+mx^{m+2}-(2+m)x^{m}+2$.\\
If $Y=0$, we will obtain $a>0.955$ for $x\in [0.4,1)$, hence $Y\neq 0$ for $a\in [0.5,0.95],\ x\in [0.4,1)$.
When $a=x=0.8$, $Y>0$, the lemma follows.\\
{\bf Lemma 3.11.} If $\rho_{1}\geq 1$, we have\\
$\prod\limits_{j=1}^{8}|\zeta_{j}|\leq (\prod\limits_{j=1}^{8}\rho_{j})(a^{2}-1+\frac{1}{4}\sum\limits_{k=1}^{8}\frac{|z_{k}|^{2}-a^{2}}{r_{k}^{2}})^{4}. $\\
Proof. By Lemma 3.2, we obtain $\sum\limits_{j=1}^{8}\frac{a^{2}-|\zeta_{j}|^{2}}{\rho_{j}^{2}} =8+2\sum\limits_{k=1}^{8}\frac{a^{2}-|z_{k}|^{2}}{r_{k}^{2}}$, then\\
$\sum\limits_{j=1}^{8}\frac{|\zeta_{j}|^{2}}{\rho_{j}^{2}}\leq 8a^{2}-8+2\sum\limits_{k=1}^{8}\frac{|z_{k}|^{2}-a^{2}}{r_{k}^{2}}.$\\
Apply the arithmetic-geometric means inequality:\\
 $\prod\limits_{j=1}^{8}\frac{|\zeta_{j}|}{\rho_{j}}
=(\prod\limits_{j=1}^{8}\frac{|\zeta_{j}|^{2}}{\rho_{j}^{2}})^{\frac{1}{2}}\leq
(\frac{1}{8}\sum\limits_{j=1}^{8}\frac{|\zeta_{j}|^{2}}{\rho_{j}^{2}})^{4}$,\\
the lemma follows.\\
We will use the following conditions
\bi
\frac{x^{2}-1}{(1-x^{m})(a+x)^{2}}+(1-a^{2})(\sigma-4)\leq 0,\ \ x\in[0.4,1].
\ei
{\bf Lemma 3.12.} If $\sigma> 4,\ \rho_{1}>1$,  and condition (3.6) holds with $m=\frac{1}{4}$, then\\
$$   4^{4}(8-\frac{9}{8}\sigma)(1-a^{2})^{-3}\leq (\sigma-4)^{4}9\prod\limits_{j=1}^{8}\rho_{j}.$$
Proof. By Lemma 3.9 and Lemma 3.11,\\
$ (8-8a^{2}-\frac{9}{8}(1-a^{2})\sigma)\prod\limits_{k=1}^{8}|z_{k}|\leq
9(\prod\limits_{j=1}^{8}\rho_{j})(a^{2}-1+\frac{1}{4}\sum\limits_{k=1}^{8}\frac{|z_{k}|^{2}-a^{2}}{r_{k}^{2}})^{4}$.\\
Write
$$\Phi=(a^{2}-1+\frac{1}{4}\sum\limits_{k=1}^{8}\frac{|z_{k}|^{2}-a^{2}}{r_{k}^{2}})(\prod\limits_{k=1}^{8}|z_{k}|)^{-\frac{1}{4}},$$
then
\bi
8-8a^{2}-\frac{9}{8}(1-a^{2})\sigma\leq 9(\prod\limits_{j=1}^{8}\rho_{j})\Phi^{4}.
\ei
Write $x=|z_{k}|,\ z_{k}=xe^{i\theta},\ r_{k}^{2}=a^{2}+x^{2}-2ax\cos \theta,\ \alpha=a^{2}-1+\frac{1}{4}\sum\limits_{l\neq k}\frac{|z_{l}|^{2}-a^{2}}{r_{l}^{2}}$,\\
we will show
\bi
(\alpha+\frac{x^{2}-a^{2}}{4r_{k}^{2}})x^{-\frac{1}{4}}\leq\alpha+\frac{1-a^{2}}{4r_{k}^{2}},
\ei
that is
$$ \frac{x^{2}-x^{m}}{1-x^{m}}-a^{2}+4\alpha r_{k}^{2}\leq 0.$$
Since \\
$\alpha\leq a^{2}-1+\frac{1}{4}\sum\limits_{l\neq k}\frac{1-a^{2}}{r_{l}^{2}}
\leq a^{2}-1+\frac{1}{4}(1-a^{2})\sigma-\frac{1-a^{2}}{4r_{k}^{2}},$\\
and $r_{k}\leq a+x$, it is sufficient to have the condition (3.6)
$$\frac{x^{2}-1}{(1-x^{m})(a+x)^{2}}+(1-a^{2})(\sigma-4)\leq 0,$$
for $\sigma>4$, we obtain (3.8).\\
Apply (3.8) eight times we have
$$ \Phi\leq \frac{1}{4}(1-a^{2})(\sigma-4).$$
Using this inequality in (3.7), the lemma follows.\\
Write $R_{k}=r_{k}\prod\limits_{j=1}^{8}r_{j}^{-\frac{1}{8}},\ $
\bi
U^{\ast}(a)=(8-v^{\ast})(\frac{2\sin\frac{\pi}{9}}{1+a})^{-\frac{7}{4}}
+(v^{\ast}-1)(\frac{1+a}{9^{\frac{1}{8}}})^{-2}+\{(\frac{2\sin\frac{\pi}{9}}{1+a})^{\frac{7}{8}(8-v^{\ast})}
(\frac{1+a}{9^{\frac{1}{8}}})^{v^{\ast}-1} \}^{2},
\ei
where
\bi v^{\ast}=\min\{j\in \mathbb{Z}: j\geq\{ 7\log(\frac{1+a}{2\sin\frac{\pi}{9}})\}
(\log\{\frac{(1+a)^{\frac{15}{8}}}{9^{\frac{1}{8}}(2\sin\frac{\pi}{9})^{\frac{7}{8}}}\})^{-1}\}.
\ei
{\bf Lemma 3.13.} If $\rho_{1}>1$,  then
$$ \sum\limits_{k=1}^{8}\frac{1}{R_{k}^{2}} \leq  U^{\ast}(a).$$
Proof. We have $ \prod\limits_{k=1}^{8}R_{k}=1$.  By (1.1), (1.2), and Lemma 3.2, we deduce
$$(\frac{2\sin\frac{\pi}{9}}{1+a})^{\frac{7}{8}}\leq R_{k}\leq \frac{1+a}{9^{\frac{1}{8}}}.$$
Taking $N=8,\ C=1,\ m=(\frac{2\sin\frac{\pi}{9}}{1+a})^{\frac{7}{8}},\
M=\frac{1+a}{9^{\frac{1}{8}}},\ c_{k}=R_{k}$ in Lemma 3.1, the lemma follows.\\
By (1.1), (1.2), Lemma 3.1 and Lemma 3.2, we obtain
\bi
\sigma\leq U(a),
\ei
where
\bi
U(a)=(8-v)(2\sin\frac{\pi}{9})^{-2}
+(v-1)(1+a)^{-2}+\{(2\sin\frac{\pi}{9})^{8-v}(1+a)^{v-1}9^{-1} \}^{2},
\ei
and
\bi
 v=\min\{j\in \mathbb{Z}: j\geq\log\{9(2\sin\frac{\pi}{9})^{-8}\}(\log\{(1+a)(2\sin\frac{\pi}{9})^{-1}\})^{-1}\}.
\ei
{\bf Lemma 3.14.} If $\rho_{1}>1,\  4<\sigma\leq U(a)<\frac{64}{9},\ a\in [0.5,0.95]$,  then
$$\frac{4U(a)}{U(a)-4}(\frac{8-\frac{9}{8}U(a)}{(1-a^{2})^{3}})^{\frac{1}{4}}\leq U^{\ast}(a).$$
Proof. Let $a\in [0.5,0.95]$, by Lemma 3.10, the following inequality is sufficient for (3.6)\\
$\lim\limits_{x\rightarrow 1-0}\frac{x^{2}-1}{(1-x^{m})(a+x)^{2}}+(1-a^{2})(\sigma-4)\leq 0,$\\
that is  $\sigma\leq 4+\frac{8}{(1-a)(1+a)^{3}},$  this is true for $\sigma<\frac{64}{9}$.\\
By Lemma 3.12, Lemma 3.13 and Lemma 3.2,\\
$ 4^{4}(8-\frac{9}{8}\sigma)(1-a^{2})^{-3}\leq \{(9\prod\limits_{j=1}^{8}\rho_{j})^{\frac{1}{4}}\sigma-4(9\prod\limits_{j=1}^{8}\rho_{j})^{\frac{1}{4}}\}^{4}\\
=\{\sum\limits_{k=1}^{8}\frac{1}{R_{k}^{2}}-4(9\prod\limits_{j=1}^{8}\rho_{j})^{\frac{1}{4}}\}^{4}\leq
 \{U^{\ast}(a)-4(9\prod\limits_{j=1}^{8}\rho_{j})^{\frac{1}{4}}\}^{4}$.\\
Using Lemma 3.12 again, we deduce
$$(8-\frac{9}{8}\sigma)^{\frac{1}{4}}(1-a^{2})^{-\frac{3}{4}}\leq \frac{1}{4}U^{\ast}(a)-4(8-\frac{9}{8}\sigma)^{\frac{1}{4}}(1-a^{2})^{-\frac{3}{4}}(\sigma-4)^{-1},$$
using (3.11), (3.12) and (3.13), the lemma follows.\\
Using (3.10) and (3.13), we obtain\\
$v^{\ast}=6: a\geq 9^{\frac{3}{17}}(2\sin\frac{\pi}{9})^{-\frac{7}{17}}-1=0.723\cdots=a_{2};$\\
$v^{\ast}=7: a< a_{2},\  a\geq (\frac{9}{2\sin\frac{\pi}{9}})^{\frac{1}{7}}-1=0.445\cdots .$\\
$v=5: a\geq 9^{\frac{1}{5}}(2\sin\frac{\pi}{9})^{\frac{5-8}{5}}-1=0.948\cdots;$\\
$v=6: a<0.948\cdots,\  a\geq 9^{\frac{1}{6}}(2\sin\frac{\pi}{9})^{\frac{6-8}{6}}-1=0.636\cdots=a_{1};$\\
$v=7: a<a_{1},\ a\geq 0.445\cdots.$\\
We divide the range of $a$ into intervals: $[0.5195, a_{1}), \ [a_{1}, a_{2}),\ [a_{2}, 0.845)$.\\
If $\rho_{1}>1$, by (3.9) and (3.12), we get contradiction with Lemma 3.14 .\\
Hence, by Lemma 3.7, we obtain Theorem 2.
\section{ Proof of the Theorem 3 }
\setcounter{equation}{0}
As in [12], for any function $h(z)$, write\\
$M(h)=\exp(\frac{1}{2\pi}\int\limits_{0}^{2\pi}\log|h(e^{i\theta})|d\theta)$.\\
Given two polynomials
$Q_{n}(z)=\sum\limits_{k=0}^{n}\left (
\begin{array}{c}
n\\
k
\end{array}
\right )a_{k}z^{k}$ and
$R_{n}(z)=\sum\limits_{k=0}^{n}\left (
\begin{array}{c}
n\\
k
\end{array}
\right )b_{k}z^{k}$,\\
we define the Szeg\"{o}-composition polynomial
$Q_{n}\otimes R_{n}=\sum\limits_{k=0}^{n}\left (
\begin{array}{c}
n\\
k
\end{array}
\right )a_{k}b_{k}z^{k}$.\\
{\bf Bruijn-Springer Theorem.} If $a_{n}b_{n}\neq 0$, we have
\bi
M(Q_{n}\otimes R_{n})\leq M(Q_{n})M(R_{n}).
\ei
Proof. See Theorem 7 of [6] and [12].\\
We will use the following well-known formula
\bi
\frac{1}{2\pi}\int\limits_{0}^{2\pi}\log|r-e^{i\theta}|d\theta=\max\{ 0, \ \log |r|\},\  \ r\in \mathbb{C}.
\ei
Let $p(z)$ be the polynomial in the introduction, write
$$p(z+a)=zQ(z+a)=\sum\limits_{k=1}^{n}\frac{Q^{(k-1)}(a)}{(k-1)!}z^{k},$$
$$p^{\prime}(z+a)=\sum\limits_{k=0}^{n-1}\frac{Q^{(k)}(a)}{k!}(k+1)z^{k},$$
$$p(z+a)=z\prod\limits_{k=1}^{n-1}(z-\gamma_{k}), \ \gamma_{k}=z_{k}-a,\ r_{k}=|\gamma_{k}|,$$
$$p^{\prime}(z+a)=n\prod\limits_{j=1}^{n-1}(z-\beta_{j}),\ \ \beta_{j}=\zeta_{j}-a,\ \rho_{j}=|\beta_{j}|.$$
We have\\
$ Q(z\rho+a)=\sum\limits_{k=0}^{n-1}\frac{Q^{(k)}(a)}{k!}\rho^{k}\left (
\begin{array}{c}
n-1\\
k
\end{array}
\right )^{-1}\left (
\begin{array}{c}
n-1\\
k
\end{array}
\right )z^{k}\\
=\sum\limits_{k=0}^{n-1}\frac{Q^{(k)}(a)}{k!}(k+1)\rho^{mk}\left (
\begin{array}{c}
n-1\\
k
\end{array}
\right )^{-1}\left (
\begin{array}{c}
n-1\\
k
\end{array}
\right )z^{k}\otimes\sum\limits_{k=0}^{n-1}\frac{1}{k+1}\rho^{(1-m)k}\left (
\begin{array}{c}
n-1\\
k
\end{array}
\right )z^{k},$ \\
hence
\bi
Q(z\rho+a)=p^{\prime}(z\rho^{m}+a)\otimes\sum\limits_{k=0}^{n-1}\frac{1}{k+1}\left (
\begin{array}{c}
n-1\\
k
\end{array}
\right )(z\rho^{1-m})^{k}.
\ei
{\bf Lemma 4.1.} We have
$$
M(p^{\prime}(z\rho+a))=n\rho^{n-1}\prod\limits_{\rho_{k}\geq\rho}\rho_{k}\rho^{-1}.
$$
Proof. By (4.2) and  $p^{\prime}(z\rho+a))=n\prod\limits_{k=1}^{n-1}(z\rho-\beta_{k}),$  we deduce \\
 $M(p^{\prime}(z\rho+a))=\exp(\frac{1}{2\pi}\int\limits_{0}^{2\pi}\{\log n+\sum\limits_{k=1}^{n-1}\log|\rho e^{i\theta}-\beta_{k}|\}d\theta)\\
 =n\rho^{n-1}\exp(\sum\limits_{k=1}^{n-1}\frac{1}{2\pi}\int\limits_{0}^{2\pi}\log|\rho^{-1}\beta_{k}-e^{i\theta}|d\theta)
 =n\rho^{n-1}\prod\limits_{\rho_{k}\geq\rho}\rho_{k}\rho^{-1},$\\
 the lemma follows.\\
{\bf Lemma 4.2.} We have
$$M( \sum\limits_{k=0}^{n-1}\frac{1}{k+1}\left (
\begin{array}{c}
n-1\\
k
\end{array}
\right )(z\rho^{1-m})^{k})=\frac{1}{n}\rho^{(1-m)(n-1)}\prod\limits_{\rho^{m-1}2\sin\frac{\pi k}{n}\geq 1}\rho^{m-1}2\sin\frac{\pi k}{n}. $$
Proof. Since\\
$\sum\limits_{k=0}^{n-1}\frac{1}{k+1}\left (
\begin{array}{c}
n-1\\
k
\end{array}
\right )(z\rho^{1-m})^{k}
=\frac{1}{n}\prod\limits_{k=1}^{n-1}(z\rho^{1-m}-(e^{2\pi i\frac{k}{n}}-1)), $\\
by (4.2), we obtain
$$ M(\sum\limits_{k=0}^{n-1}\frac{1}{k+1}\left (
\begin{array}{c}
n-1\\
k
\end{array}
\right )(z\rho^{1-m})^{k})=\frac{1}{n}\rho^{(1-m)(n-1)}\prod\limits_{\rho^{m-1}2\sin\frac{\pi k}{n}\geq 1}\rho^{m-1}2\sin\frac{\pi k}{n},$$
the lemma follows.\\
But
\bi
Q(z\rho+a)=\prod\limits_{k=1}^{n-1}(z\rho-\gamma_{k}),\ M(Q(z\rho+a))=\rho^{n-1}\prod\limits_{r_{k}\geq\rho}r_{k}\rho^{-1},
\ei
 by Lemma 4.1, Lemma 4.2, (4.1), (4.3) and (4.4), Theorem 3 follows.\\
 {\bf Remark 4.1.} It is possible to obtain new results on the Sendov conjecture by combining Theorem 3 with Lemma 3.1.

{ {\small E-mail:\ mengzzh@126.com}
\end{document}